\documentclass[12pt]{amsart}

\usepackage{amssymb}

\newtheorem{Thm}{Theorem}[section]

\newtheorem{Lem}[Thm]{Lemma}
\newtheorem{Prop}[Thm]{Proposition}

\theoremstyle{definition}

\newtheorem{Rem}[Thm]{Remark}

\theoremstyle{remark}

\newcommand{\vfi}{\varphi}
\def\la{\lambda}

\def \eps{\varepsilon}

\def\Ndb{\mathbb N}

\def\Rdb{\mathbb R}

\begin{document}
\title{Embeddings of locally finite metric spaces into Banach spaces}

\author {F. Baudier}
\address{Universit\'e de Franche-Comt\'e, Laboratoire de Math\'ematiques UMR 6623,
16 route de Gray, 25030 Besan\c con Cedex, FRANCE.}
\email{florent.baudier@math.univ-fcomte.fr}

\author {G. Lancien}
\address{Universit\'e de Franche-Comt\'e, Laboratoire de Math\'ematiques UMR 6623,
16 route de Gray, 25030 Besan\c con Cedex, FRANCE.}
\email{glancien@math.univ-fcomte.fr}

\begin{abstract} We show that if $X$ is a Banach space without cotype, then
every locally finite metric space embeds metrically into $X$.

\end{abstract}

\subjclass[2000]{46B20 (primary), 51F99 (secondary)}

\maketitle

\section{Introduction}\label{intro}

Let $X$ and $Y$ be two Banach spaces. If $X$ and $Y$ are linearly isomorphic,
the {\it Banach-Mazur distance} between $X$ and $Y$, denoted by $d_{BM}(X,Y)$,
is the infimum of $\|T\|\,\|T^{-1}\|$, over all linear isomorphisms $T$ from
$X$ onto $Y$.

For $p\in [1,\infty]$ and $n\in \Ndb$, $\ell_p^n$ denotes the space $\Rdb^n$
equipped with the $\ell_p$ norm. We say that a Banach space $X$ uniformly
contains the $\ell_p^n$'s if there is a constant $C \geq 1$ such that for every
integer $n$, $X$ admits an $n$-dimensional subspace $Y$ so that
$d_{BM}(\ell_p^n,Y)\leq C$.

\smallskip A metric space $M$ is {\it locally finite} if any ball of $M$ with
finite radius is finite. If moreover, there is a function $C:(0,+\infty)\to
\Ndb$ such that any ball of radius $r$ contains at most $C(r)$ points, we say
that $M$ has a {\it bounded geometry}.

Let $(M,d)$ and $(N,\delta)$ be two metric spaces and $f:M\to N$ be a map. For
$t>0$ define:
$$\rho_f(t)=\inf\{\delta(f(x),f(y)),\ \ d(x,y)\geq t\}$$
and
$$\omega_f(t)=\sup\{\delta(f(x),f(y)),\ \ d(x,y)\leq t\}.$$
We say that $f$ is a {\it coarse embedding} if $\omega_f(t)$ is finite for all
$t>0$ and $\lim_{t\to \infty}\rho_f(t)=\infty$.

\noindent Suppose now that $f$ is injective. We say that $f$ is a {\it uniform
embedding} if $f$ and $f^{-1}$ are uniformly continuous. Following \cite{MN},
we also define the {\it distortion} of $f$ to be
$$ {\rm dist}(f):= \|f\|_{Lip}\|f^{-1}\|_{Lip}=\sup_{x\neq y \in
M}\frac{\delta(f(x),f(y))}{d(x,y)}.\sup_{x\neq y \in
M}\frac{d(x,y)}{\delta(f(x),f(y))}.$$ If the distortion of $f$ is finite, we
say that $f$ is a {\it metric embedding} and that $M$ metrically embeds into
$N$ and we denote $M {\hookrightarrow} N$. If dist$(f)\leq C$, we use the
notation $M \buildrel {C}\over {\hookrightarrow} N$.

\medskip
A celebrated result of I. Aharoni \cite{A} asserts that every separable metric
space metrically embeds into $c_0$ (the space of all real sequences converging
to $0$, equipped with the supremum norm). It is an open problem to know what
are the Banach spaces that share this property with $c_0$. In other words, if
$c_0$ metrically embeds into a separable Banach space $X$, do we have that $X$
contains a closed subspace which is linearly isomorphic to $c_0$? Very
recently, N.J. Kalton \cite{K} made an important step in this direction by
showing that one of the iterated duals of such a Banach space has to be non
separable (in particular, $X$ cannot be reflexive). In fact, he even proved it
for $X$ such that $c_0$ coarsely or uniformly embeds into $X$.

In another direction, N. Brown and E. Guentner proved in \cite{BG} that any
metric space with bounded geometry coarsely embeds into a reflexive Banach
space. This was also improved by Kalton in \cite{K}, who obtained that any
locally finite metric space uniformly and coarsely embeds into a reflexive
space. These results are connected with the coarse Novikov conjecture. Indeed,
Kasparov and Yu have shown in \cite{KY} that if a metric space with bounded
geometry coarsely embeds into a super-reflexive Banach space, then it satisfies
this important conjecture. In this paper, we improve the results of Brown,
Guentner and Kalton, by showing that any Banach space uniformly containing the
$\ell_\infty^n$'s is metrically universal for all locally finite metric spaces.
It should also be noted, that Brown and Guentner used a specific space,
$(\sum\ell_{p_n})_2$, with $p_n$ tending to infinity and that this space
fulfils our hypothesis.

\section{Results}\label{results}

\begin{Thm}\label{main} There exists a universal constant $C>1$ such that
for every Banach space $X$ uniformly containing the $\ell_\infty^n$'s and every
locally finite metric space $(M,d)$: $M \buildrel {C}\over {\hookrightarrow}
X$.
\end{Thm}

\begin{proof} Let $X$ be a Banach space uniformly containing the
$\ell_\infty^n$'s (or equivalently without any non trivial cotype). Our first
lemma follows directly from the classical work of B. Maurey and G. Pisier
\cite{MP}.

\begin{Lem}\label{fcod} For any finite codimensional subspace $Y$ of $X$, any $\eps>0$ and
any $n\in \Ndb$, there exists a subspace $F$ of $Y$ such that
$d_{BM}(\ell_\infty^n,F)<1+\eps$.
\end{Lem}

We shall also need the following version of Mazur's Lemma (see for instance
\cite{LT1} page 4, for a proof).

\begin{Lem}\label{mazur} Let $X$ be an infinite dimensional Banach space, $E$ be a
finite dimensional subspace of $X$ and $\eps>0$. Then there is a finite
codimensional subspace $Y$ of $X$ such that:
$$\forall x\in E\ \ \forall y\in Y,\ \ \|x\|\leq (1+\eps)\|x+y\|.$$
\end{Lem}

Consider now a locally finite metric space $(M,d)$. For $t\in M$ and $r\geq 0$,
$B(t,r)=\{s\in M,\ d(s,t)\leq r\}$. We fix a point $t_0$ in $M$ and denote
$|t|=d(t,t_0)$, for $t\in M$. Since $M$ is locally finite, we may assume, by
multiplying $d$ by a constant if necessary, that $B(t_0,1)=\{t_0\}$. For any
non negative integer $n$, we denote $B_n=B(t_0,2^{n+1})$. Using Lemmas
\ref{fcod} and \ref{mazur}, together with the fact that each ball $B_n$ is
finite, we can build inductively finite dimensional subspaces
$(F_n)_{n=0}^\infty$ of $X$ and $(T_n)_{n=0}^\infty$ so that for every $n\geq
0$, $T_n$ is a linear isomorphism from $\ell_\infty(B_n)$ onto $F_n$ satisfying
$$\forall u\in \ell_\infty(B_n)\ \ \ \ \frac{1}{2}\|u\|\leq \|T_nu\|\leq \|u\|$$
and also such that $(F_n)_{n=0}^\infty$ is a Schauder finite dimensional
decomposition of its closed linear span $Z$. More precisely, if $P_n$ is the
projection from $Z$ onto $F_0\oplus...\oplus F_n$ with kernel $\overline{\rm
sp}\,(\bigcup_{i=n+1}^\infty F_i)$, we will assume as we may, that $\|P_n\|\leq
2$. We denote now $\Pi_0=P_0$ and $\Pi_n=P_n-P_{n-1}$ for $n\geq 1$. We have
that $\|\Pi_n\|\leq 4$.

We now consider $\varphi_n:B_n \to \ell_\infty(B_n)$ defined by
$$\forall t\in B_n,\ \ \varphi_n(t)= \big(d(s,t)-d(s,t_0)\big)_{s\in B_n}=
\big(d(s,t)-|s|\big)_{s\in B_n}.$$ The map $\varphi_n$ is clearly an isometric
embedding of $B_n$ into $\ell_\infty(B_n)$.

\noindent Then we set:
$$\forall t\in B_n,\ \ f_n(t)=T_n(\varphi_n(t)) \in F_n.$$

\noindent Finally we construct a map $f:M\to X$ as follows:

(i) $f(t_0)=0$.

(ii) For $n\geq 0$ and $2^n\leq |t|<2^{n+1}$:
$$f(t)=\lambda
f_n(t)+(1-\lambda)f_{n+1}(t),\ \ {\rm where}\ \
\lambda=\frac{2^{n+1}-|t|}{2^{n}}.$$

\noindent We will show that dist$(f)\leq 9\times 24=216$.

\noindent Note first that for any $t$ in $M$, $$\frac{1}{16}|t|\leq
\|f(t)\|\leq |t|.$$

\medskip We start by showing that $f$ is Lipschitz. Let $t,t' \in M\setminus \{t_0\}$ and assume,
as we may, that $1\leq |t|\leq |t'|$ (we recall that $B(t_0,1)=\{t_0\}$).

\medskip\noindent I) If $|t|\leq \frac{1}{2}|t'|$, then
$$\|f(t)-f(t')\|\leq |t|+|t'|\leq \frac{3}{2}|t'|\leq 3(|t'|-|t|)\leq
3\,d(t,t').$$

\medskip\noindent II) If $\frac{1}{2}|t'|<|t|\leq |t'|$, we have two different
cases to consider.

\smallskip 1) $2^n\leq |t|\leq |t'|<2^{n+1}$, for some $n\geq 0$. Then,
let
$$\lambda=\frac{2^{n+1}-|t|}{2^{n}}\ \ {\rm and}\ \ \lambda'=\frac{2^{n+1}-|t'|}{2^{n}}.$$
We have that
$$|\lambda-\lambda'|=\frac{|t'|-|t|}{2^{n}}\leq \frac{d(t,t')}{2^{n}},$$
so

\begin{equation*}
\begin{split}
\|f(t)-f(t')\|&=\|\lambda f_n(t)-\lambda'
f_n(t')+(1-\lambda)f_{n+1}(t)-(1-\lambda')f_{n+1}(t')\|\\
&\le \lambda\|f_n(t)-f_n(t')\|+(1-\lambda)\|f_{n+1}(t)-f_{n+1}(t')\|+2|\lambda-\lambda'|\,|t'|\\
&\le d(t,t')+2^{n+2}|\lambda-\lambda'|\le 5\,d(t,t').
\end{split}
\end{equation*}

\smallskip 2) $2^n\leq |t|<2^{n+1}\leq |t'|<2^{n+2}$, for some $n\geq 0$. Then,
let
$$\lambda=\frac{2^{n+1}-|t|}{2^{n}}\ \ {\rm and}\ \ \lambda'=\frac{2^{n+2}-|t'|}{2^{n+1}}.$$
We have that
$$\lambda \leq \frac{d(t,t')}{2^n},\ \ {\rm so}\ \ \lambda|t|\leq 2\,d(t,t').$$
Similarly
$$1-\lambda'=\frac{|t'|-2^{n+1}}{2^{n+1}}\leq \frac{d(t,t')}{2^{n+1}}\ \ {\rm
and}\ \ (1-\lambda')|t'|\leq 2\,d(t,t').$$

\begin{equation*}
\begin{split}
\|f(t)&-f(t')\|=\|\lambda f_n(t)+(1-\lambda)f_{n+1}(t)-\lambda'
f_{n+1}(t')-(1-\lambda')f_{n+2}(t')\|\\
&\le
\lambda(\|f_n(t)\|+\|f_{n+1}(t)\|)+(1-\lambda')(\|f_{n+1}(t')\|+\|f_{n+2}(t')\|)+\|f_{n+1}(t)-f_{n+1}(t')\|\\
&\le d(t,t')+2\lambda|t|+2(1-\lambda')|t'|\le 9\,d(t,t').
\end{split}
\end{equation*}

\smallskip\noindent We have shown that $f$ is 9-Lipschitz.

\medskip We shall now prove that $f^{-1}$ is Lipschitz. We consider
$t,t'\in M\setminus\{t_0\}$ and assume again that $1\leq |t|\leq |t'|$. We need
to study three different cases. In our discussion, whenever $|t|$ (respectively
$|t'|$) will belong to $[2^{m},2^{m+1})$, for some integer $m$, we shall denote
$$\lambda=\frac{2^{m+1}-|t|}{2^{m}}\ \ ({\rm respectively}\ \ \lambda'=\frac{2^{m+1}-|t'|}{2^{m}}).$$

\smallskip 1) $2^{n}\leq |t|\leq |t'|<2^{n+1}$, for some $n\geq 0$. Then

$$\Pi_n(f(t)-f(t'))=T_n(\lambda\vfi_n(t)-\lambda'\vfi_n(t'))\ \ {\rm and}$$
$$[\lambda\vfi_n(t)-\lambda'\vfi_n(t')]_{|_{s=t'}}=\la d(t,t')+(\la' -\la)|t'|.$$
So
$$2\|\Pi_n(f(t)-f(t'))\|\geq \la d(t,t')+(\la' -\la)|t'|.$$
On the other hand
$$\Pi_{n+1}(f(t)-f(t'))=T_{n+1}((1-\lambda)\vfi_{n+1}(t)-(1-\lambda')\vfi_{n+1}(t'))\ \ {\rm and}$$
$$[(1-\lambda)\vfi_{n+1}(t)-(1-\lambda')\vfi_{n+1}(t')]_{|_{s=t'}}=(1-\lambda) d(t,t')+(\la -\la')|t'|.$$
So
$$2\|\Pi_{n+1}(f(t)-f(t'))\|\geq (1-\lambda) d(t,t')+(\la -\la')|t'|.$$
Therefore
$$16\|f(t)-f(t')\|\geq d(t,t').$$

\smallskip 2) $2^n\leq |t|<2^{n+1}\leq |t'|<2^{n+2}$, for some $n\geq 0$.

$$2\|\Pi_{n}(f(t)-f(t'))\|=2\la\|T_{n}(\vfi_{n}(t))\|\geq \la|t|.$$

$$2\|\Pi_{n+2}(f(t)-f(t'))\|=2(1-\la')\|T_{n+2}(\vfi_{n+2}(t'))\|\geq (1-\la')|t'|.$$

$$\Pi_{n+1}(f(t)-f(t'))=T_{n+1}((1-\lambda)\vfi_{n+1}(t)-\la'\vfi_{n+1}(t'))\ \
{\rm and}$$
$$[\la'\vfi_{n+1}(t')-(1-\lambda)\vfi_{n+1}(t)]_{|_{s=t}}=\la' d(t,t')-\la'|t|+(1-\la)|t|.$$
So
$$2\|\Pi_{n+1}(f(t)-f(t'))\|\geq \la' d(t,t')-\la'|t|+(1-\la)|t|.$$
Combining our three estimates, we obtain
$$24\|f(t)-f(t')\|\geq \la' d(t,t')+(1-\la')(|t|+|t'|)\geq d(t,t').$$

\smallskip 3) $2^n\leq|t|<2^{n+1}<2^p\leq |t'|<2^{p+1}$ for some integers $n$
and $p$.

$$2\|\Pi_{p}(f(t)-f(t'))\|=2\la'\|T_{p}(\vfi_{p}(t'))\|\geq \la'|t'|\ \ {\rm and}$$

$$2\|\Pi_{p+1}(f(t)-f(t'))\|=2(1-\la')\|T_{p+1}(\vfi_{p+1}(t'))\|\geq (1-\la')|t'|.$$
So
$$24\|f(t)-f(t')\|\geq \frac{3}{2}|t'|\geq |t'|+|t|\geq d(t,t').$$

\medskip\noindent All possible cases are settled and we have shown that $f^{-1}$ is 24-Lipschitz.

\end{proof}

We conclude by mentioning the following easy fact.

\begin{Prop} Let $N$ be a metric space. The following assertions
are equivalent:

(i) For any locally finite metric space $M$: $M {\hookrightarrow} N$.

(ii) There exists $C\geq 1$ such that for any locally finite metric space $M
\buildrel {C}\over {\hookrightarrow} N$.

\end{Prop}

This is an immediate consequence of the following lemma.

\begin{Lem} Let $(M_p,d_p)_{p=1}^\infty$ be a sequence of locally finite metric
spaces. Then there exists a locally finite metric space $(M,d)$ such that each
$M_p$ embeds isometrically into $M$.
\end{Lem}

\begin{proof} For any $p\in \Ndb$, pick $x_0^p \in M_p$. Consider
$M=\bigcup_{p=1}^\infty \{p\}\times M_p$. Let $x\in M_p$ and $y\in M_q$. We
define $d((p,x),(q,y))=d_p(x,y)$ if $p=q$ and
$d((p,x),(q,y))=\max\{p,q,d_p(x_0^p,x),d(x_0^q,y)\}$ if $p\neq q$. We leave it
to the reader to check that $(M,d)$ is a locally finite metric space.

\end{proof}

\begin{Rem} We do not know if the converse of Theorem \ref{main} is true. So
the question is: if every locally finite metric space metrically embeds in a
given Banach space $X$, do we have that $X$ uniformly contains the
$\ell_\infty^n$'s? However, it follows from the work of Mendel and Naor in
\cite{MN} that such a Banach space cannot be K-convex, or equivalently, it must
contain the $\ell_1^n$'s uniformly.
\end{Rem}


\begin{thebibliography}{WW}

\bibitem{A} I. Aharoni, Every separable metric space is Lipschitz equivalent to
a subset of $c_0^+$, {\it Israel J. Math.}, {\bf 19} 1974, 284-291.

\bibitem{BG} N. Brown and E. Guentner, Uniform embeddings of bounded geometry
spaces into reflexive Banach spaces, {\it Proc. Amer. Math. Soc.}, {\bf 133(7)}
2005, 2045-2050.

\bibitem{K} N.J. Kalton, Coarse and uniform embeddings into reflexive spaces,
preprint.

\bibitem{KY} G. Kasparov and G. Yu, The coarse Novikov conjecture and uniform
convexity, {\it Advances Math.}, to appear.

\bibitem{LT1} J. Lindenstrauss and L. Tzafriri, Classical Banach spaces I,
Springer Berlin 1977.

\bibitem{MP} B. Maurey, G. Pisier, S\'{e}ries de variables al\'{e}atoires vectorielles
ind\'{e}pendantes et propri\'{e}t\'{e}s g\'{e}om\'{e}triques des espaces de Banach, {\it Studia
Math.}, {\bf 58(1)} 1976, 45-90.

\bibitem{MN} M. Mendel, A. Naor, Metric cotype, preprint.



\end{thebibliography}
\end{document}